\documentstyle[12pt]{report}

\def\R{I\kern-.30em{R}}
\def\N{I\kern-.30em{N}}
\def\F{I\kern-.30em{F}}
\def\C{I\kern-.60em{C}}
\def\Z{Z\kern-.50em{Z}}

\def\P{I\kern-.30em{P}}
\def\E{I\kern-.30em{E}}

\def\build#1_#2^#3{\mathrel{\mathop{\kern 0pt#1}\limits_{#2}^{#3}}}

\def\Sum{\displaystyle\sum}

\def\supp{\mbox{\rm supp}\ }
\def\dist{\mbox{\rm dist}\ }

\def\nth{n$^{\mbox{\footnotesize th}}$ }

\def\ess{\mbox{\footnotesize ess}}

\newtheorem{theorem}{Theorem}[section]
\newtheorem{lemma}{Lemma}[section]
\newtheorem{proposition}{Proposition}[section]
\newtheorem{corollary}{Corollary}[section]
\newtheorem{definition}{Definition}[section]

\def\eop{{\vrule height5pt width5pt depth0pt}\par\bigskip}

\newenvironment{prf}{\noindent{\it Proof}}{\eop}
\newenvironment{prfoft}{\noindent{\it Proof of Theorem}}{\eop}

\newcommand{\rf}[1]{(\ref{#1})}
\newcommand{\beq}[1]{\begin{equation}\label{#1}}
\newcommand{\eeq}{\end{equation}}
\newcommand{\ba}{\begin{array}}
\newcommand{\ea}{\end{array}}

\font\twelve=cmbx10 at 15pt
\font\ten=cmbx10 at 12pt

\begin{document}

\begin{titlepage}

\begin{center}

\renewcommand{\thefootnote}{\fnsymbol{footnote}}

{\ten Centre de Physique Th\'eorique\footnote{Unit\'e Propre de
Recherche 7061} - CNRS - Luminy, Case 907} {\ten F-13288 Marseille
Cedex 9 - France }

\vspace{0.7 cm}

{\twelve LANDAU HAMILTONIANS \\
WITH RANDOM POTENTIALS : LOCALIZATION \\ AND THE DENSITY OF STATES}

\vspace{0.3 cm}

\setcounter{footnote}{0}
\renewcommand{\thefootnote}{\arabic{footnote}}

{\bf J.M. Combes\footnote{
D\'epartement de Math\'ematiques, Universit\'e de Toulon et du Var
- 83130 La Garde - France.}$^,$\footnote{
Supported in part by CNRS.} and P.D. Hislop\footnote{ Mathematics
Department, University of Kentucky - Lexington KY 40506-0027 -
USA.}$^,$\footnote{
Supported in part by NSF grants INT 90-15895 and DMS 93-07438.} }

\vspace{0.7 cm}

{\bf Abstract}

\end{center}

We prove the existence of localized states at the edges of the
bands for the two-dimensional Landau Hamiltonian with a random
potential, of arbitrary disorder, provided that the magnetic field
is sufficiently large. The corresponding eigenfunctions decay
exponentially with the magnetic field and distance. We also prove
that the integrated density
of states is Lipschitz continuous away from the Landau energies.
The proof relies on a Wegner estimate for the finite-area magnetic
Hamiltonians with random potentials and exponential decay estimates
for the finite-area Green's functions. The proof of the decay
estimates for the Green's functions uses
fundamental results from two-dimensional bond percolation theory.

\vspace{0.7 cm}

\noindent Key-Words : Landau Hamiltonians, random operators,
localization.

\bigskip

\noindent Number of figures : 4

\bigskip

\noindent August 1994

\noindent CPT-94/P.3061

\bigskip

\noindent anonymous ftp or gopher: cpt.univ-mrs.fr

\end{titlepage}


\renewcommand{\thechapter}{\arabic{chapter}}
\renewcommand{\thesection}{\thechapter}

\setcounter{chapter}{1}
\setcounter{equation}{0}

\section{Introduction}\label{S.1}

The existence of localized states for a two-dimensional gas of
non-interacting electrons in a constant magnetic field is a main
ingredient in various discussions and proofs of the integer quantum
Hall effect (see e.g.~\cite{[Halp]},
\cite{[Prange-Joynt]}, \cite{[Kunz]}, \cite{[Bellissard]} and
\cite{[Prange-Girvin]}). It is generally believed that localization
occurs near the band edges for large magnetic fields and bounded,
random potentials of arbitrary disorder. According to Halperin's
argument~\cite{[Halp]}, the localization length should diverge near
the Landau levels. This is in contrast to the situation with no
magnetic field. For two dimensional random systems, localization is
expected to hold at all energies for arbitrary disorder and the
eigenfunctions are expected to decay exponentially.

In this paper, we study the family $ H_{ \omega } $ of
two-dimensional Landau
Hamiltonians with Anderson-type potentials, having mean zero, on
$L^2 (\R^2)$.
We prove that localization does occur in all energy intervals
$I_{n}( B ) \equiv [ (2n + 1) B + {\cal O} (B^{-1}),\ (2n + 3) B -
{\cal O} (B^{-1}) ], n = 0, 1, 2, \ldots $ at large magnetic field
strengths B and for arbitrary disorder. Recall that $\sigma (H_{
\omega })$ is contained in bands about the Landau levels $E_{n}( B
) \equiv (2n + 1) B,\ n = 0,1,2,\dots$, of width $|| V
||_{\infty}$, independent of B. We follow the approach
of~\cite{[CoHi]} developed to study random Schr\"odinger operators
on $L^2 (\R^d)$. This work~\cite{[CoHi]} extends to continuous
systems the techniques of Howland~\cite{[H]},
Simon-Wolff~\cite{[SiWo]}, von Dreifus-Klein~\cite{[VDk]}, and
Spencer~\cite{[S]}.

For large magnetic fields, we justify a one Landau band
approximation (of~\cite{[Weg]}) uniformly in $ n $ and obtain
exponential decay estimates in $x$ and $B$ on the Green's function
for finite-area Hamiltonians. The key to these estimates is showing
that equipotential lines don't percolate with high probability. For
potentials with zero averages, this holds at all energies except
the Landau levels, which correspond to the critical percolation
threshold. In addition to this restriction, there is a small region
of energy of ${\cal O} (B^{-1})$ around each Landau level where
small denominators in the interband perturbation expansion can't be
controlled by our method. Although this is in agreement with an
earlier conjecture of Laughlin~\cite{[L]}, it remains an open
question whether these small bands of energies about the Landau
levels correspond to extended states.

In section~\ref{S.2} below, we describe the model and state the
main results. We also give some elementary estimates needed later
to justify the one Landau band approximation. In section~\ref{S.3},
we prove Wegner estimates for the quantum Hall Hamiltonian
restricted to finite boxes. As a by-product, we obtain the
Lipschitz continuity of the integrated density of states away from
the Landau energies. The proofs of the exponential decay of the
finite-area Green's function are given in section~\ref{S.4}. The
results of section~\ref{S.3} on the Wegner estimate and
section~\ref{S.4} on the decay of the finite volume Hamiltonians
are used in section~\ref{S.5}, together with results
of~\cite{[CoHi]}, to prove the main theorem. We prove some
technical lemmas in the appendix.

We have recently learned of some related results on localization
for the models studied here by J.\ Pul\'{e} ~\cite{[Pu]} and by W.\
M.\ Wang ~\cite{[W2]}.

\section*{Acknowledgments}

We thank T.~Hoffmann-Ostenhof and W.~Thirring for their hospitality
at the Erwin Schr\"odinger Institute in Vienna where most of this
work was done. We thank R.~Seiler and V.~Jaksic for helpful
discussions on the quantum Hall effect.


\setcounter{chapter}{2}
\setcounter{equation}{0}

\section{The Model and the Main Results}\label{S.2}

We consider a one-particle Hamiltonian which describes an electron
in two-dimensions $(x_1, x_2)$ subject to a constant magnetic field
of strength $B > 0$ in the perpendicular $x_3$-direction, and a
random potential $V_{\omega}$. The Hamiltonian $H_{\omega}$ has the
form \beq{2.1}
H_{\omega} = (p - A)^2 + V_{\omega},
\eeq

\noindent on the Hilbert space $L^2 (\R^2)$, where $p \equiv - i
\nabla$, and the vector potential $A$ is \beq{2.2}
A = {B \over 2} (x_2, -x_1),
\eeq

\noindent so the magnetic field $B = \nabla \times A$ is in the
$x_3$-direction. The random potential $V_{\omega}$ is Anderson-like
having the form
\beq{2.3}
V_{\omega} (x) = \Sum_{i \in \Z^d} \lambda_i (\omega) u (x - i).
\eeq

\noindent We make the following assumptions on the single-site
potential $u$ and the coupling constants $\{ \lambda_i (\omega)
\}$.

\begin{itemize}

\item[(V1)] $u\ge 0,\ u \in C^2,\ \supp u \subset B (0, {1 \over
\sqrt 2})$, and $\exists C_0 > 0$ and $r_0 > 0$ s.t. $u | B (0,
r_0) > C_0$.

\item[(V2)] $\{ \lambda_i (\omega) \}$ is an independent,
identically distributed family of random variables with common
distribution $g \in C^2 ([- M, M ])$ for some $0 < M < \infty$,
s.t. $\int \lambda g (\lambda) d \lambda = 0$ and $g (\lambda) > 0$
Lebesgue a.e. $\lambda \not= 0$.

\end{itemize}

We denote by $H_A \equiv (p - A)^2$, the Landau Hamiltonian. As is
well-known, the spectrum of $H_A$ consists of an increasing
sequence $\{ E_n (B) \}$ of eigenvalues, each of infinite
multiplicity, given by
\beq{2.4}
E_n (B) = (2 n + 1) B,\ n = 0,1,2,\dots
\eeq

\noindent Note that $D (H_{\omega}) = D (H_A)\ \forall\ \omega \in
\Omega$. We will call $E_n (B)$ the \nth Landau level and denote by
$P_n$ the projection onto the corresponding subspace. The
orthogonal projection is denoted by $ Q_{n} \equiv 1 - P_{n} $. Let
$M_0 \equiv\ \build{\sup}_{x, \omega}^{} || V_{\omega} || <
\infty$. Then, $\sigma (H_{\omega}) \subset \build{\cup}_{n \ge
0}^{} \sigma_n$, where $\sigma_n \equiv [ E_n (B) - M_0,\ E_n (B) +
M_0 ]$, which we call the \nth Landau band. We show that
$\sigma (H_{\omega})$ is deterministic.
The magnetic translations are
defined for $a \in \Z^2$ by
\beq{2.5}
U_a \equiv e^{-i B x \wedge a} e^{-i p \cdot a} \eeq

\noindent where $ x \wedge a \equiv x_2 a_1 - x_1 a_2$. We then
have \beq{2.6}
U_a H_{\omega} U_a^{-1} = H_{T_a \omega}, \eeq

\noindent where $T_a : \Omega \to \Omega$ is the
$\Z^2$-translation. Standard results (cf.~\cite{[CL]}) show that
$H_{\omega}$ is a $\Z^2$-ergodic self-adjoint family of operators
and consequently its spectrum is deterministic. Note that
$\sigma (H_{\omega})$ is not necessarily equal a.s. to
$\build{\cup}_{n \ge 0}^{} \sigma_n$. Provided some $V_{\omega},\
\omega \in \Omega_0$ and $| \Omega_0 | > 0$, lifts the degeneracy
of the Landau level, then ergodicity implies that the spectrum
consists of bands each of which lies in some
interval about $E_n (B)$, which might be strictly contained in
$\sigma_n$.

\begin{theorem}\label{T.2.1}

Let $H_{\omega}$ be the family given in~\rf{2.1} with vector
potential
$A$ satisfying~\rf{2.2}, $B > 0$, and the random potential
$V_{\omega}$ as in~\rf{2.3} and satisfying~(V1)-(V2). Let $$
I_n (B) \equiv \left[ E_n (B) + {\cal O} (B^{-1}),\ E_{n+1} (B) -
{\cal O} (B^{-1}) \right].
$$

\noindent There exists $B_0 \gg 0$ such that for $B > B_0$ and all
$n = 0,1,2,\dots,$
$$
\sigma (H_{\omega}) \cap I_n (B)
$$
is pure point almost surely and the corresponding eigenfunctions
decay exponentially. The integrated density of states is Lipschitz
continuous away from $\sigma (H_A)$.

\end{theorem}

Let us make two remarks about the theorem. First, we note that the
above theorem holds at arbitrary disorder. For large disorder, the
techniques of ~\cite{[CoHi]} apply directly to show that, without
the percolation estimates, $ \sigma ( H_{ \omega } ) $ is almost
surely pure point in each Landau band. This regime, however, is of
little interest as the quantum Hall conductivity vanishes in this
case.
Secondly, we show, in fact, that the localization length, for
energies near the band edges as in Theorem 2.1, is a decreasing
function of the field strength $ B $ so that the wave functions are
strongly localized. We also show that the localization length
increases as the energy approaches the Landau levels. The precise
manner in which this occurs follows from Proposition 5.1 and
Theorem 5.2. However, our method fails to give an estimate of the
power law divergence of the localization length near the Landau
level.

As is clear from the Wegner estimate, Theorem~\ref{T.3.1}, our
method fails to give information about the integrated density of
states at the Landau energies. However, we can improve the result
if we make a
stronger hypothesis in~(V1) on $\supp u$.

\begin{corollary}\label{C.2.2}

If, in addition to the hypothesis of Theorem~\ref{T.2.1}, we have
$u \ge C_0 \chi_{\Lambda_1 (0)},\ C_0 > 0$, then the integrated
density of states is Lipschitz continuous.

\end{corollary}

If the hypothesis of Corollary 2.1 does not hold, then a large
portion of configuration space is unaffected by the potential. It
is not, therefore, surprising that there is a discontinuity in the
integrated density of states at the Landau energies as there is
for the Landau Hamiltonian. A phenomenon of this type has been
observed by Brezin et. al. ~\cite{[BZJ]} for a Poisson distribution
of impurities at low energy. So we do not expect that the IDS is
Lipschitz continuous at the Landau energies without a condition of
the support of $ u $ which implies that the zero set of $ V_{
\omega } $ is in some sense ``small''.

We mention that W.M.~Wang~\cite{[Wang]} has obtained an asymptotic
expansion in the semi-classical limit for the density of states at
large magnetic field strengths away from the Landau levels,
partially justifying the one-band approximation.

We conclude this section with some simple observations on the
Landau projections $P_n$.

The projection $P_n$ on the \nth Landau level of $H_A$ has a kernel
given by
\beq{2.7}
P_{n} (x,y) = B e^{- {i B \over 2} x \wedge y} p_n \left( B^{1
\over 2} (x - y) \right),
\eeq

\noindent where $p_n (x)$ is of the form \beq{2.8}
p_n (x) = \left\{ \mbox{\nth degree polynomial in}\ x \right\}\
e^{- {| x |^2 \over 2} },
\eeq

\noindent and independent of $B$. We will make repeated use of the
following elementary lemma, the proof of which follows by direct
calculation using the kernel~\rf{2.7}-\rf{2.8}.

\begin{lemma}\label{L.2.3}

Let $\chi_1,\chi_2$ be functions of disjoint, compact support with
$| \chi_i | \le 1$ and let $\delta \equiv \dist (\supp \chi_1,\
\supp \chi_2) > 0$. Then,

\begin{itemize}

\item[(1)] $|| \chi_1 P_n \chi_1 ||_1 \le C_n B | \supp \chi_1 |^2
\ ;$

\item[(2)] $|| \chi_1 P_n \chi_1 ||_{H S} \le C_n B | \supp \chi_1
|\ ;$

\item[(3)] $|| \chi_1 P_n \chi_2 ||_{H S} \le C_n B | p_n \left(
B^{1 \over 2} \delta \right) |^{1 \over 2} \left\{ | \supp \chi_1 |
| \supp \chi_2 | \right\}^{1 \over 2} \ ;$

\item[(4)] $|| \chi_1 P_n ||_{H S} \le C_n B,$

\end{itemize}

\noindent where $C_n$ varies from line to line and depends only on
$\chi_i$ and $n$, and $ H S $ denotes the Hilbert Schmidt norm.
Note that $| p_n \left( B^{1 \over 2} \delta \right) |^{1
\over 2} \le C_0 e^{- \epsilon B}$, for any $\epsilon > 0$ and $B$
large enough.

\end{lemma}


\setcounter{chapter}{3}
\setcounter{equation}{0}

\section{Wegner Estimate}\label{S.3}

We define local Hamiltonians as relatively compact perturbations of
the Landau Hamiltonian $H_A = (p - A)^2$, as defined in section 2.
Let $\Lambda
\subset \R^2$ denote an open connected region in $ \R^2 $. We let
$\Lambda_{\ell} (x)$ denote a square of side $\ell$ centered at $x
\in \R^2$,
$$
\Lambda_{\ell} (x) \equiv \left\{ y \in \R^2 |\ | x_i - y_i | <
\ell,\ i = 1,2 \right\}.
$$

\noindent Given $\Lambda \subset \R^2$, the local potential
$V_{\Lambda}$ is defined as follows. Freeze all $\lambda_j (\omega)
\in \Z^2 \cap (\R^2 \backslash \Lambda)$ and consider $\widetilde
V$ so obtained. This potential depends on the external, fixed
coupling constants and on all $\lambda_i (\omega) \in \Z^2 \cap
\Lambda$. We define
$V_{\Lambda} \equiv \widetilde V | \Lambda$ and define $H_{\Lambda}
\equiv H_A + V_{\Lambda}$ on $L^2 (\R^2)$. These Hamiltonians are
not independent of the external configurations but we will prove
estimates uniform in the external random variables. We will use the
conditional probability law
\beq{3.1}
\overline P (A \cap B) \le \overline P (A) \overline P (B), \eeq

\noindent where $\overline P$ is the probability conditioned on the
external variables and $A\ \&\ B$ are any two events in $\Lambda$.
Note that $\sigma_{\ess}(H_{\Lambda}) = \sigma_{\ess} (H_A)$, since
$V_{\Lambda}$ is relatively compact.

We prove the following theorem.

\begin{theorem}\label{T.3.1}

$\exists B_0 > 0$ and a constant $C_W > 0$ such that for all $B >
B_0$ and for any $E \not\in \sigma (H_A)$, \[ \P \{ \dist (\sigma
(H_{\Lambda}),\ E) < \delta \} \le C_W [ \dist (\sigma
(H_{\Lambda}),\ E) ]^{- 2} || g ||_{\infty} \delta B | \Lambda |.\]

\end{theorem}

This theorem will follow from the properties of the spectral
projectors for $H_A$ and a spectral averaging theorem. Since
$H_{\Lambda}$ depends analytically on the coupling constants
$\lambda_i$, we need only a simple version which we state without
proof (cf.~\cite{[CoHi]}, \cite{[KS]}, \cite{[CHM]}).

\begin{lemma}\label{L.3.2}

Let $H_{\lambda} \equiv H_0 + \lambda u,\ \lambda \in \R$, be a
self-adjoint family on $D (H_0)$ with $C_0 D^2 \le u \le M <
\infty$. Let $E_{\lambda} (\cdot)$ be the spectral family for
$H_{\lambda}$. For any $h \ge 0$, $\supp h$ compact, and $h \in
L^{\infty} (\R)$, and for any $L \subset \R$ measurable, we have $$
|| \int_{\R} h (\lambda) D E_{\lambda} (L) D d \lambda || \le
C_0^{-1} | L |\ || h ||_{\infty}.
$$

\end{lemma}

For simplicity, we will work with the case $ n = 0 $, the first
Landau band, although
the calculation is uniform in $ n $. As can be easily checked, the
calculations depend only on the difference between the energy $ E $
that we are considering and the nearest Landau energy $ E_{ n }( B
) $ and hence is independent of $ n $. To begin the proof of
Theorem~\ref{T.3.1}, we need a simple estimate. Let $ \Delta $ be
an interval in the first Landau band $ \sigma_0 $. Let $ E_{ \Delta
} $ be the spectral projector for $ H_{\Lambda} $ associated with $
\Delta $.

\begin{lemma}\label{L.3.3}

$|| E_{\Delta} Q_0 E_{\Delta} || \le d_{\Delta}^{- {1 \over 2} }
\left( 1 - (2 d_{\Delta})^{-1} | \Delta | \right)^{- {1 \over 2} }
M^{1 \over 2}$, where \hfill \break
$d_{\Delta} \equiv \dist (\sigma (H_A) \backslash \{ B \},\ \Delta)
= {\cal O} (B)$.

\end{lemma}

\begin{prf}

Let $E_m \in \Delta$ be the center of the interval. We then can
write $$
\ba{ll}
E_{\Delta} Q_0 E_{\Delta} & \le \dist (\sigma (H_A) \backslash \{ B
\},\ \Delta)^{-1} (E_{\Delta} (H_A - E_m) Q_0 E_{\Delta}) \\[3mm] &
\le d_{\Delta}^{-1} \{ E_{\Delta} (H_{\Lambda} - E_m) Q_0
E_{\Delta} + E_{\Delta} V_{\Lambda} Q_0 E_{\Delta} \}. \ea
$$
This implies that
$$
|| E_{\Delta} Q_0 E_{\Delta} || \le d_{\Delta}^{-1} \left\{ {|
\Delta | \over 2} || E_{\Delta} Q_0 E_{\Delta} || + M_{0} || Q_0
E_{\Delta} || \right\}.
$$

\noindent Since $d_{\Delta} = {\cal O} (B)$, it is clear that for
all $B$ sufficiently large $(2 d_{\Delta})^{-1} | \Delta | \ll 1$,
so $$
|| E_{\Delta} Q_0 E_{\Delta} || \le d_{\Delta}^{-1} \left( 1 - (2
d_{\Delta})^{-1} | \Delta | \right)^{-1} M_{0} || E_{\Delta} Q_0
E_{\Delta} ||^{1 \over 2},
$$
and the result follows. \end{prf}

Note that as $d_{\Delta} = {\cal O} (B)$, we obtain \beq{3.2}
|| E_{\Delta} Q_0 E_{\Delta} || = {\cal O} \left( B^{- {1 \over 2}
} \right).
\eeq

\begin{prfoft}~\ref{T.3.1}

We can assume without less of generality that the closest point in
$\sigma (H_A)$ to $E$ is $E_0 (B) = B$. All the calculations below
hold for any band. Let $\Delta \subset \sigma_0 \backslash \{ E_0
(B) \}$ be a connected interval containing $E$ and let $E_{\Delta}$
be the spectral projection for $H_{\Lambda}$ and $\Delta$. Recall
{}from Chebyshev's inequality that
\beq{3.3}
\P_{\Lambda} \{ \dist (\sigma (H_{\Lambda}),\ E) < \eta \} \le
\E_{\Lambda} (Tr E_{\Delta}),
\eeq

\noindent where $\P_{\Lambda}$ and $\E_{\Lambda}$ denote the
probability and expectation with respect to the variables in
$\Lambda \cap \Z^2$ and
$Tr$ denotes the trace on $L^2 (\R^2)$. We first note that
\beq{3.4}
Tr E_{\Delta} \le 2 Tr (P_0 E_{\Delta} P_0). \eeq

\noindent This follows from the identity $$
Tr E_{\Delta} = Tr E_{\Delta} P_0 E_{\Delta} + Tr E_{\Delta} Q_0
E_{\Delta},
$$
and the bound
$$
Tr E_{\Delta} Q_0 E_{\Delta} \le || E_{\Delta} Q_0 E_{\Delta} ||
(Tr E_{\Delta}),
$$

\noindent since $E_{\Delta} Q_0 E_{\Delta} \ge 0$. Now by Lemma
3.2, $|| E_{\Delta} Q_0 E_{\Delta} || = {\cal O} \left( B^{- {1
\over 2} } \right)$ so~\rf{3.4} follows for all $B$ sufficiently
large. Let us now suppose $\inf \Delta > B$ for definiteness.
{}From~\rf{3.4}, and positivity we obtain \beq{3.5}
\ba{ll}
Tr E_{\Delta} P_0 E_{\Delta} & \le Tr ( E_{\Delta} (H_{\Lambda} -
B) P_0 (H_{\Lambda} - B) E_{\Delta} ) \cdot \dist (\Delta, B)^{-2}
\\[3mm]
& \le Tr ( P_0 V_{\Lambda} E_{\Delta} V_{\Lambda} P_0 ) \cdot \dist
(\Delta, B)^{-2}.
\ea
\eeq

\noindent Writing $V_{\Lambda} = \Sum_i \lambda_i u_i$ for short,
the trace in~\rf{3.5} is
\beq{3.6}
\Sum_{i, j} \lambda_i \lambda_j Tr (P_0 u_i E_{\Delta} u_j P_0),
\eeq

\noindent where $i,j \in \Lambda \cap \Z^2$. Defining $A^{i j}
\equiv u_i^{1 \over 2} A u_j^{1 \over 2}$ for any $A \in B ({\cal
H})$, we have from~\rf{3.6},
\beq{3.7}
\Sum_{i, j} \lambda_i \lambda_j Tr \left( P_0^{i j} E_{\Delta}^{i
j} \right).
\eeq
We must estimate
\beq{3.8}
\ba{ll}
\E_{\Lambda} \left( \Sum_{i, j} \lambda_i \lambda_j Tr \left(
P_0^{i j} E_{\Delta}^{i j} \right) \right) & \le \Sum_{i, j}
\E_{\Lambda} \left( | \lambda_i \lambda_j |\ | Tr \left( P_0^{i j}
E_{\Delta}^{i j}
\right) | \right) \\[3mm]
& \le {1 \over 2} M^2 \Sum_{i, j} || P_0^{i j} ||_1 \E_{\Lambda}
\left\{ || E_{\Delta}^{i i} || + || E_{\Delta}^{j j} || \right\}.
\ea
\eeq

\noindent Since $E_{\Delta}^{i i} \equiv u_i^{1 \over 2} E_{\Delta}
u_i^{1 \over 2}
\ge 0$ and self-adjoint, we have
$$
\ba{ll}
\E_{\Lambda} \left( || E_{\Delta}^{i i} || \right) & \le
\build{\sup}_{\psi,\ || \psi || = 1}^{} \left\{ \E_{\Lambda} \left(
\langle \psi,\ E_{\Delta}^{i i} \psi \rangle \right) \right\}
\\[3mm] & \le C_0^{-1} || g ||_{\infty} | \Delta |, \ea
$$
by Lemma~\ref{L.3.2}. Consequently, \rf{3.8} is bounded above by
\beq{3.9}
{1 \over 2} C_0^{-1} M^2 || g ||_{\infty} | \Delta | \Sum_{i, j} ||
P_0^{i j} ||_1.
\eeq

\noindent To evaluate the trace norm, we first note that for $i =
j,\ P_0^{i i} > 0$ so by Lemma~\ref{L.2.3}, \beq{3.10}
\Sum_{i = j} || P_0^{i j} ||_1 = C_0 B | \supp u |^2 | \Lambda |.
\eeq

\noindent Next, suppose $u_i u_j \not= 0,\ i \not= j$. Let $\chi_{i
j}$ be the characteristic function on $\supp u_i u_j$. Then, if $i
\cap j$ denotes the set of such pairs,
\beq{3.11}
\Sum_{i \cap j} || u_i^{1 \over 2} P_0 u_j^{1 \over 2} ||_1 \le
\Sum_{i
\cap j} || \chi_{i j} P_0 \chi_{i j} ||_1 \le C_1 B |\Lambda|\ |
\supp u |^2,
\eeq

\noindent where we used $\sup u_i \le 1$. Finally, for $u_i u_j =
0$, let $\{
\chi_{\ell}^2 \}$ be a partition of unity covering $\Lambda$ so
that $\chi_{\ell} | \supp u_{\ell} = 1$ and $\chi_{\ell} \chi_n =
0,\ \ell \not= m$. Using the inequality
\beq{3.12}
|| A B ||_1 \le ||A ||_{H S} || B ||_{H S}, \eeq
we obtain
\beq{3.13}
|| u_j^{1 \over 2} P_0 u_j^{1 \over 2} || \le \Sum_{\ell} || u_j^{1
\over 2} P_0
\chi_{\ell} ||_{H S} || \chi_{\ell} P_0 u_j^{1 \over 2} ||_{H S}.
\eeq

\noindent As in Lemma~\ref{L.2.3} with $\delta \equiv | i - \ell |
- 2 r_u$, we easily compute
\beq{3.14}
|| u_i^{1 \over 2} P_0 \chi_{\ell} ||_{H S} \le C_2 B e^{- B ( | i
-
\ell | - 2 r_u )},
\eeq

\noindent where $0 < r_u < {1 \over \sqrt 2}$ is the radius of
$\supp u_i$ (and $\chi_{\ell}$) (see~(V1)). Summing over $\{ i j
\}'$, the set of pairs with $u_i u_j = 0$, we get
{}from~\rf{3.13}-\rf{3.14}, \beq{3.15}
\ba{ll}
\Sum_{ \{ i j \}' } || P_0^{i j} ||_1 & \le C_3 B^2 \Sum_{ \{ i j
\}' \ell } e^{- B ( | i - \ell | - 2 r_u )} e^{- B ( | j - \ell | -
2 r_u)} \\[3mm]
& \le C_4 B^2 | \Lambda | e^{- \delta B}, \ea
\eeq

\noindent for same $\delta > 0$. Combining~\rf{3.10}, \rf{3.11} and
\rf{3.15} in \rf{3.9} we obtain an upper bound for all $B$ large
enough,
$$
\E (Tr E_{\Delta}) \le C_W || g ||_{\infty} | \Delta | B | \Lambda
|, $$

\noindent where $C_W$ depends on $M^2,\ C_0^{-1}$, and $|
\supp u |^2$. This proves the theorem. \end{prfoft}

The estimate of Theorem~\ref{T.3.1} suffices to prove the Lipschitz
continuity of the integrated density of states away from the Landau
levels, as stated in Theorem~\ref{T.2.1}. With regard to
Corollary~\ref{C.2.2}, let us show how the additional hypothesis on
$\supp u$ allows the improvement. For $M_{0} \equiv || V_{ \omega }
||_{ \infty } $ as in section~\ref{S.2}, define
$$
H_0 = H_A + 2 M_{0} (1 - \chi_{\Lambda}), $$
and the finite-area Hamiltonian by
$$
H_{\Lambda} = H_0 + V_{\Lambda}.
$$

\noindent Beginning with~\rf{3.4}, we have for $\Delta \subset
\sigma_0 $
and $E_m
\equiv$ center of $\Delta$,
$$
\ba{ll}
Tr E_{\Delta} & \le 2 Tr \left\{ E_{\Delta} (H_A + 2 M_{0} - E_m)
P_0 (B + 2 M_{0} - E_m)^{-1} \right\} \\[3mm]
& \le 2 (B + 2 M_{0} - E_m)^{-1} \left\{ Tr E_{\Delta} (H_{\Lambda}
- E_m) P_0 \right. \\[3mm]
& \left. \hphantom{\le} + Tr E_{\Delta} (2 M_{0} \chi_{\Lambda} -
V_{\Lambda}) P_0
\right\}.
\ea
$$

\noindent Since $2 M_{0} \chi_{\Lambda} - V_{\Lambda} > M_{0}
\chi_{\Lambda}$ and $|| (H_{\Lambda} - E_m) E_{\Delta} || \le {|
\Delta | \over 2}$, we obtain
$$
Tr E_{\Delta} \le 2 (B + 2 M_{0} - E_m)^{-1} \left\{ {| \Delta |
\over 2} Tr P_0 E_{\Delta} + M_{0} Tr E_{\Delta} \chi_{\Lambda} P_0
\right\}. $$
As $(B + 2 M_{0} - E_m)^{-1} | \Delta | < {| \Delta | \over M_{0}}
\ll 1$, we arrive at
$$
Tr E_{\Delta} \le 4 M_{0} C_1^{-1} (B + 2 M_{0} - E_m)^{-1} \left\{
\Sum_{i \in \Lambda \cap \Z^2} Tr (E_{\Delta} u_i P_0) \right\}. $$

\noindent Here we used the fact that $\Sum_{i \in \Z^2 \cap
\Lambda} u_i \ge C_1 \chi_{\Lambda}$. The remaining steps are the
same as above. In light of this calculation one might speculate
that the singularity at the Landau energies of the IDS is due to
the existence of large regions where there is no potential. Indeed,
numerical studies on the Poisson model~\cite{[BZJ]} seem to also
support this idea.


\setcounter{chapter}{4}
\setcounter{equation}{0}

\section{Percolation Theory and Decay Estimates}\label{S.4}

In this section, we prove the technical estimates required to
justify the one-Landau band approximation. We consider for
simplicity the first Landau band $\sigma_0 \equiv [ B - M_{0},\ B +
M_{0} ]$, but
all other bands can be analysed using the same techniques. The
results are uniform in the band index $n$. Formally, if one
neglects the band interaction, the effective Hamiltonian for an
electron at energy $E$ is $E = B + V (x)$. Consequently, in this
approximation, the electron motion is along equipotential lines $V
(x) + B - E = 0$. Since $V$ is random, it is natural to estimate
the probability that these equipotential lines percolate through a
given box. If not, the electron will remain confined to bounded
regions. One can expect that the interband interaction will not
change this picture. We will do this in the second part of this
section by showing that the Green's function decays exponentially
in $x$ and $B$ through regions where $| V (x) + B - E | > a > 0$.
The first part of this section is devoted to reformulating our
problem as a problem in bond percolation.

\subsection{Percolation Estimates}

We first show that in annular regions between boxes of side $\ell$
and $\ell / 3$ there exist closed, connected ribbons where the
condition $| V (x) + B - E | > a > 0$ is satisfied, provided $E\neq
B$, with a probability which converges exponentially fast to 1 as
$\ell$ tends to infinity. Obviously, the existence of such a ribbon
is equivalent to the impossibility for equipotential lines at
energy $E$ to percolate from the center of the box to its boundary.
Although this is a classical matter, let us recall how one can
formulate the above condition in terms of two-dimensional bond
percolation.

Recall that $V_{\omega} (x) = \Sum_{i \in \Z^2} \lambda_i (\omega)
u (x - i)$, where the single-site potential $u \ge 0$ and has
support inside a ball of radius $r_u < {1 \over \sqrt 2}$. We
define $ r_u $ to be the smallest radius such that $ supp~u \subset
B(0,r_u) $. Consider a new
square lattice $\Gamma \equiv e^{i \pi / 4} \sqrt 2 \Z^2$. The
midpoint of each bond of $\Gamma$ is a site of $\Z^2$ (see
Figure~1). We will denote by $b_j$ the bond of $\Gamma$ having $j
\in \Z^2$ as it's midpoint. For definiteness, we assume $E \in (B,\
B + M_0)$. The other energy interval can be treated similarly.

\begin{definition}\label{D.4.1}

The bond $b_j$ of ~$\Gamma$ is {\bf occupied} if $\lambda_j
(\omega) < {E - B \over 2}$. The probability $\P \left\{
\vphantom{{E - B \over 2}} \lambda_j (\omega) \right.$ $\left. < {E
- B \over 2} \right\} \equiv p$ is the probability that $b_j$ is
occupied ($p$ is independent of $j$ by the $i i d$ assumption).

\end{definition}

Let us assume that the bond $b_j$ is occupied and consider (see
Figure~2),
\beq{4.1}
{\cal R}_j \equiv \left\{ x | \dist (x,\ b_j) < {1 \over \sqrt 2} -
r_u \equiv r_1 \right\}.
\eeq

\noindent Obviously, ${\cal R}_j$ does not intersect the support of
the other single-site potentials centered on $\Z^2 \backslash \{ j
\}$ so that $V (x) = \lambda_j (\omega) u (x - j)\ \forall \ x \in
{\cal R}_j$. Then, if $b_j$ is occupied, one has $V (x) < {E - B
\over 2}
\ \forall \ x \in {\cal R}_j$ (recall that ${E - B \over 2} > 0$).
We now assume that there is a closed circuit of occupied bonds
${\cal C} \equiv \build{\cup}_{j \in \gamma}^{} b_j,\ \gamma
\subset \Z^2$ (i.e. a connected union occupied bonds). We call $ {
\cal R } \equiv \build{\cup}_{j \in \gamma}^{} {\cal R}_j $ the
closed {\em ribbon }
associated with ${\cal C}$. For all $x \in \cal{R}$, we have $V (x)
< {E - B \over 2}$. If we take $a \equiv {E - B \over 2}$, then
\beq{4.2}
V (x) + B - E < - a\ \forall \ x \in \cal{R}. \eeq

The existence of a closed ribbon $\cal{R}$ so that $V$ satisfies
condition~\rf{4.2} is a consequence of the existence of a closed
circuit ${\cal C}$ in $\Gamma$ of occupied bonds. In order to
estimate the probability that ${\cal C}$ exists, we use some
standard results of percolation theory (see, e.g.~\cite{[Gr]} and
\cite{[ChCh]}) which we now summarize.

Let $\Z^2$ be the square lattice (the length of the side plays no
role in the calculations). A {\it bond} (edge) of $\Z^2$ is said to
be {\it occupied} with probability $p,\ 0 \le p \le 1$, and {\it
empty} with
probability $1 - p$. We are interested in the case when the bonds
are independent (called Bernoulli bond percolation). The critical
percolation probability $p_c$ is defined as follows. Let
$P_{\infty} (p)$ be the probability that the origin belongs to an
infinite (connected) cluster of occupied bonds. Then, we define $$
p_c \equiv \inf \{ p | P_{\infty} (p) > 0 \}. $$

\noindent For 2-dimensional Bernoulli bond percolation, $p_c = {1
\over 2}$. Hence if $p > p_c$, occupied bonds percolate ; that is,
we can find a connected cluster of occupied bonds running off to
infinity with non-zero probability.

Of importance for us are the results concerning the existence of
closed circuits of occupied bonds. Let $r_{n,\ell}$ be a rectangle
in $\Z^2$ of width $\ell$ and length $n \ell$. Let $R_{n,\ell}$ be
the probability that there is a crossing of $r_{n,\ell}$, the long
way, by a connected path of occupied bonds. A basic result is

\begin{theorem}\label{T.4.2}

For $p > p_c,\ R_{n,\ell} \ge 1 - C_0 n \ell e^{- m (1 - p) \ell}$,
for some constant $C_0$.

\end{theorem}

The exponential factor $m (q)$ is strictly positive for $q < p_c$
and

$m (q) \searrow 0$ as $q \nearrow p_c$. This factor measures the
probability that the origin $0$ is connected to $x \in \Z^2$ by a
path of occupied bonds
$$
P_{0 x} (p) \le e^{- m (p) |x|}.
$$

Let us write $r_\ell$ for $r_{1,\ell}$, the box of side $\ell$. An
annular region between two concentric boxes is denoted by $a_{\ell}
\equiv r_{3 \ell} \backslash r_\ell$. A closed circuit of occupied
bonds in $a_\ell$ is a connected path of occupied bonds lying
entirely within $a_\ell$. Using Theorem~\ref{4.1} and the $F G K$
inequality, one can compute the probability $A_\ell$ of a closed
circuit of occupied bonds in $a_\ell$ for $p > p_c$.

\begin{theorem}\label{T.4.3}

For any $p \in [0,1],\ A_\ell \ge [ R_{3,\ell} (p)]^4$. In
particular, if $p > p_c,\ \exists 0 < C_0 < \infty$ as in Theorem
4.1, such that \beq{4.3}
A_{\ell} \ge 1 - 12 C_0 e^{- m (1 - p) \ell}. \eeq

\end{theorem}

We now apply these results to our situation as follows. On the
lattice $\Gamma$ defined above, the probability that any bond is
occupied is given by
$$
p = \int_{- M}^a g (\lambda) d \lambda,
$$

\noindent so, under our assumptions on the density $g$, if $a > 0$
then $p > p_c = {1 \over 2} $, and we are above the critical
percolation threshold
$p_c = {1 \over 2}$. Note that when $E = B,\ a = 0$ so $p = {1
\over 2} = p_c$, the critical probability. It follows from
Theorem~\ref{T.4.3} that any annular region $a_{\ell} \equiv r_{3
\ell} \backslash r_\ell$ in
$\Gamma$ of diameter $\sqrt 2 {\ell} \equiv {1 \over 2} (3 \sqrt 2
{ \ell} -
\sqrt 2 {\ell})$ and sides parallel to the bonds of $\Gamma$ (see
Figure~3) contains a closed circuit of occupied bonds with
probability given by~\rf{4.3}. By the argument above, there is a
ribbon $\cal{R}$ associated with ${\cal C}$ in $a_\ell$ whose
properties we summarize in the next proposition.

\begin{proposition}\label{P.4.4}

Assume (V1) and (V2). Let $r_u \equiv {1 \over 2} \mbox{diam}\
(\supp u),\ {\ell} > \sqrt 2,\ E \in \sigma_0 \backslash \{ B \}$,
and $a > 0$. Then for $m (q)$ and $C_0$ as in Theorem~\ref{T.4.3},
$\exists$ a ribbon $\cal{R}$ satisfying
\beq{4.4}
\mbox{diam}{\cal R} \ge 2 \left( {1 \over \sqrt 2} - r_u \right) ;
\eeq

\beq{4.5}
\ba{ll}
& \dist ({\cal R},\ \partial r_{3 \ell}),\ \dist ({\cal R},\
\partial r_ \ell) \ge {1
\over \sqrt 2} + r_u ; \\[3mm]
& \hphantom{\dist (} {\cal R} \subset a_\ell, \ea
\eeq
and s.t.
\beq{4.6}
V (x) + B - E < - a,\ \forall\ x \in {\cal R}, \eeq
with a probability larger than
\beq{4.7}
1 - 12 C_0 e^{- m (1 - p) \ell},
\eeq
where
\beq{4.8}
p \equiv \int_{- M}^a g (\lambda) d \lambda. \eeq

\end{proposition}

\def\doublespaced{\baselineskip=\normalbaselineskip
\multiply\baselineskip by 2}
\def\doublespaced{doublespaced}

\def\doublespaced{\baselineskip=\normalbaselineskip
\multiply\baselineskip by 2}
\def\doublespace{\doublespaced}

\def\C{\makebox[.8em][l]{\makebox[.25em][l]{C}\rule{.03em}{1.5ex}}}

\newcommand{\bt}[1]{\par\smallskip\noindent{\bf #1}\sl}
\newcommand{\et}{\par\smallskip\goodbreak}
\newcommand{\bp}[1]{\par\smallskip\noindent{\it#1}}
\newcommand{\ep}{ \hfill QED\smallskip\goodbreak \par}

\subsection{Decay Estimates}

The effective one Landau band
Hamiltonian $ B + V $ localizes electrons at energies $E$ where the
equipotential lines
$ E = V (x) + B $ don't percolate to infinity. The effect of the
interband interaction is to induce some tunneling through the
``Classically Forbidden'' ribbon $\cal{R}$ of Proposition 4.1. As a
consequence, instead of localization
in the compact subsets of $ I \hspace{-.17cm} R ^2 $ bounded by
$\cal{R}$, one expects exponential decay of the Green's function
in $x$ and $B$ across such ribbons $\cal{R}$. Such an estimate is
the starting point of the inductive, multi-scale analysis detailed
in section 5. By the geometric resolvent equation, we show there
that it suffices to consider the following ideal situation, where
for some $ a > 0 $,
\beq{4.9}
V(x) + B - E < -a , ~~ \forall \: x \in I \hspace{-.17cm} R ^ 2 ,
\eeq
or, alternately,
\beq{4.10}
V(x) + B - E > a , ~~ \forall \: x \in
I \hspace{-.17cm} R ^ 2 .
\eeq
A condition such as (4.9) with $ E > B $ is satisfied, with a
probability given in Proposition 4.1, by a smoothing (see section
5) of the potential $ V_{\cal{R}} $ defined as \beq{4.11}
V_{\cal{R}} (x) = \left \{ \begin{array}{ll} V(x) & x \in \cal{R}
\\
0 & x \in I \hspace{-.17cm} R ^2 \setminus \cal{R} . \end{array}
\right.
\eeq
Here we obtain decay estimates on
\[ H = H_A + V \]
with $V$ having compact support with non-empty interior and
satisfying (4.9) or (4.10).

Let $ {\cal O} $ be an open, bounded, connected set in $ I
\hspace{-.17cm} R ^2 $ with smooth boundary and define $ \rho (x)
= \mbox{dist} \: ( x , {\cal O} ) $. Let $ \eta \in C_0^{ \infty }
( I \hspace{-.17cm} R ^2 ) $ with $ \eta > 0 $ and $ \mbox{supp} \:
\eta \subset B_1 (0) $. For any $ \epsilon > 0 $, define $ \eta _{
\epsilon } (x) = \eta ( x / \epsilon ) $. We consider the smoothed
distance function $ \rho _{ \epsilon } (x) \equiv ( \eta _{
\epsilon }
\star \rho ) (x) $; $ \mbox{supp} \: \rho _{ \epsilon } \subset I
\hspace{-.17cm} R ^2 \setminus \{ x | \mbox{dist} \: ( x , {\cal
O}^{c} ) < \epsilon \} $. We fix $ \epsilon > 0 $ small and
write $ \rho $ for
$ \rho _{ \epsilon } $ below for simplicity. We have $ || \nabla
\rho ||_{ \infty } < C_0 / \epsilon $ and $ || \Delta \rho ||_{
\infty } < C_1 / \epsilon ^2, $ for constants $ C_0 , ~ C_1 > 0 $
depending only on $ \eta $ and $ {\cal O} $. This $ \epsilon $ will
play no role in the analysis below and, consequently, we absorb it
into the constants $ C_{0} $ and $ C_{1} $. We consider
one-parameter families of operators defined for $ \alpha \in I
\hspace{-.17cm} R $ as
\beq{4.12}
H_A ( \alpha ) \equiv e^{ i \alpha \rho } H_A e^{ - i \alpha \rho }
;
\eeq
\vspace{-.3in}
\beq{4.13}
H ( \alpha ) \equiv H_A ( \alpha ) + V ; \eeq
\vspace{-.3in}
\beq{4.14}
P ( \alpha ) \equiv e^{ i \alpha \rho } P e^{ -i \alpha \rho } , ~~
\mbox{etc.},
\eeq

\noindent
and similarly for the local Hamiltonian
$ H_{ \Lambda } ( \alpha )
\equiv H_A ( \alpha ) + V _{ \Lambda } $. Here, we write $ P $ for
the projector $ P_{0} $. For
$ \alpha \in I \hspace{-.17cm} R $, these families are unitarity
equivalent with the $ \alpha = 0 $ operator. \vspace{.1in}

\noindent
{\bf Lemma 4.1}. {\em The family $ H ( \alpha ) ~~(and ~similarly
{~}for ~H_{ \Lambda } ( \alpha ) ) , ~~ \alpha \in I \hspace{-.17cm
} R $,
has an analytic continuation into the strip \beq{4.15}
S \equiv \{ \alpha \in \C | \:
| \mbox{Im} \: \alpha | <
\eta_{\rho}B^{ 1/2 }, \}
\eeq
as a type $A$ analytic family with domain $ D (H) $. The positive
constant $ \eta_{\rho} $ depends only on the distance function $
\rho $. Furthermore, in this strip $S$, one has $ P ( \alpha ) ^2 =
P ( \alpha ) $ and for some constant $ C_1 $ independent of $
\alpha $, \beq{4.16}
|| P ( \alpha ) || < C_1
\eeq
and}
\beq{4.17}
|| Q ( \alpha ) ( H_A ( \alpha ) - z ) ^{-1} || < C_1 B ^{-1} , ~~
\mbox{\em if} ~~ \mbox{dist} \: ( z , B ) \leq B . \eeq
\vspace{.1in}

\noindent
{\bf Proof}. For $ \alpha \in I \hspace{-.17cm} R $, one has
\beq{4.18}
\begin{array}{rl}
H_A ( \alpha ) & = ( - i \nabla - \alpha \nabla \rho - A ) ^2 \\ &
\\
& = H_A - \alpha [ \nabla \rho \cdot ( p - A ) + ( p - A ) \cdot
\nabla \rho ] + \alpha ^2 | \nabla \rho | ^2 \\ & \\
& = H_A + \alpha ^2 | \nabla \rho | ^2 + i \alpha \Delta \rho - 2
\alpha \nabla \rho \cdot ( p - A ) . \end{array}
\eeq
By a standard unitary equivalence argument, it suffices to show
that
\[ \{ \alpha ^2 | \nabla \rho | ^2 + i \alpha \Delta \rho - 2
\alpha \nabla \rho \cdot ( p - A ) \}( H_A - z ) ^{-1}, \] has norm
less than 1 for some $ z \not \in \sigma ( H_A ) $ and $ \alpha $
purely imaginary (cf \cite{[K]}). For later
purposes, we choose $ z \in C (B) \equiv \{ z | \: | z - B | = B \}
$, circle of
radius $B$ centered at $B$.
Since $ | \nabla \rho | ^2 $ and $ \Delta \rho $ are bounded, we
can choose $ \eta_{\rho} ' $ such that \[ || (\alpha ^2 | \nabla
\rho | ^2 + i \alpha \Delta \rho ) ( H_A - z ) ^{-1} || < 1/2, \]
for $ | \alpha | < \eta_{\rho} ' B^{ 1/2 } $ and for all $ z \in
C(B) $. Hence,
it is enough to show that $ \forall \: | \alpha | < \eta_{ \rho
}''B^{ 1/2 } $, for a possibly smaller constant $ \eta_{ \rho }''
$, \beq{4.19}
2 | \alpha | \: || \nabla \rho \cdot
( p - A ) ( H_A - z ) ^{-1} || < 1/2 ,
\eeq
for some $ z \in C (B) $.
Since $ || \nabla \rho ||_{ \infty }
< C_0 $, we easily find that $ || \nabla \rho \cdot ( p - A ) ( H_A
- z ) ^{-1} || < C_1 B^{ - 1/2 } $. This implies (4.19) for all $ B
$ sufficiently large. We take $ \eta_{ \rho } $ to be the smallest
of these two constants. From these estimates and
(4.18) for $ \alpha \in S $ we have that for $B$ large enough,
\beq{4.20}
|| ( H_A ( \alpha ) - z ) ^{-1} || <
C_2 B ^{-1} , ~~ z \in C (B) ,
\eeq
for some constant $ C_2 $ uniform in $ \alpha \in S $ and $ z \in
C(B) $. Next, note that the eigenfunctions of $ H_A $ are analytic
vectors for the family $ e^{ i \alpha \rho } , ~ \alpha \in S $. It
is a consequence of this, the analyticity of $ H_A $,
and the eigenvalue equation, that $ \sigma ( H_A ( \alpha ) ) $ is
independent of $ \alpha , ~ \alpha \in S $. The family $ P ( \alpha
) , ~ \alpha \in I \hspace{-.17cm} R $, has an analytic
continuation in $S$ given by the contour integral \beq{4.21}
P ( \alpha ) = \frac{ -1 }{ 2 \pi i } \int_{ C(B) } ( H_A ( \alpha
) - z ) ^{-1} \: d z .
\eeq
The boundedness of $ P ( \alpha ) $ follows from (4.20) and (4.21).
The idempotent property of $ P( \alpha ) $ follows from this
analyticity and the identity $ P( \alpha )^2 = P( \alpha ) $, which
holds for real $ \alpha $.
Furthermore, the function $ \alpha \in S \rightarrow Q( \alpha )(
H_A ( \alpha ) - z ) ^{-1} $ is holomorphic on and inside $ C (B)
$. By the maximum modulus principle, it follows that
\beq{4.22}
|| Q ( \alpha ) ( H_A ( \alpha ) - z ) ^{-1} || \leq {\displaystyle
\sup_{ z \in C(B) } } \: || Q ( \alpha ) ( H_A ( \alpha ) - z )
^{-1} ||, \eeq
and the bound (4.17) follows from this, (4.16) and (4.20). $ \Box $

We next prove the main estimate of this section. \vspace{.1in}

\noindent
{\bf Theorem 4.3}. {\em Assume that $ ( V , E , B ) $ satisfy (4.9)
or (4.10) for some $ a > 0 $ and
$ E \in \sigma_0 \setminus \{ B \} $. Furthermore, assume that $
supp ~V $ is compact with non-empty interior. There exists
constants $C_2 \leq \eta_{ \rho }$, $C_3$, and $B_1$, depending
only on $ M_0 \equiv || V ||_{ \infty } , ~
|| \nabla \rho ||_{ \infty },
and ~ || \nabla V ||_{ \infty } , $
such that if we define $ \gamma
\equiv C_2 \: \mbox{min} \: \{ B^{ 1/2 } , a B \}, $ and $u$ is a
solution of
\beq{4.23}
( H_A + V - z ) u = v , ~~ z \equiv E + i \epsilon , ~~ \epsilon
00, ~~ E > 0,
\eeq
for some $ v \in D ( e^{ \gamma \rho} ) $, then for $ B > B_1 , ~~
\forall \: \alpha \in \C , ~~ | \mbox{Im} \: \alpha | < \gamma $,
we have \beq{4.24}
u \in D( e^{i \rho \alpha}),
\eeq
\beq{4.25}
|| e^{ \alpha \rho } P u || \leq
C_{3} a ^{-1} || e^{ \alpha \rho } v || , \eeq
and}
\beq{4.26}
|| e^{ \alpha \rho } Qu || \leq C_{3} B ^{-1} || e ^{ \alpha \rho }
v || .
\eeq
\vspace{.1in}

\noindent
{\bf Proof}. Let $ v( \alpha ) = e^{i \alpha \rho }v $, so that $
v( \alpha ) $ is analytic in the strip $ | Im ( \alpha )| < \gamma
$. Let
$ u( \alpha ) = e^{i \alpha \rho} u, ~~\alpha \in \R $, i.e.,
\beq{4.27}
u(\alpha) = ( H( \alpha ) - z )^{-1} v( \alpha ) \equiv ( ( H - z
)^{-1} v)( \alpha),
\eeq
\noindent with $ H \equiv H_{A} + V $, as above. Since $ V $ is $
H_{A} $-compact by assumption, $ H $ has point spectrum, and
according to Lemma 4.1 and standard arguments, $ H( \alpha ) $ has
real spectrum independent of $ \alpha $ in the strip $ S $ defined
in (4.15).
Then, $ u( \alpha ) $ has an analytic continuation in $ | Im (
\alpha )| < \gamma $, which proves (4.24).

Projecting the equation (4.23) for $u$ along $ P ( \alpha ) $ gives
\beq{4.28}
( B + V - z ) ( P u ) ( \alpha ) = ( P v ) ( \alpha ) + ( [ Q V P -
P V Q ] u ) ( \alpha ) .
\eeq
Taking the scalar product of (4.28) with $ ( P u ) ( \alpha ) $
results in the inequality, \beq{4.29}
\begin{array}{rl}
a || ( P u ) ( \alpha ) || ^2 & \leq || ( P u ) ( \alpha ) || \: ||
( P v ) ( \alpha ) || \\ & \\
& + \{ || ( P^{*} Q ) ( \alpha ) || \: || ( Q V P ) ( \alpha ) ||
\} || ( P u ) ( \alpha ) ||^2 \\ & \\ & + || ( P V Q ) ( \alpha )
|| \: || ( P u ) ( \alpha ) || \: || ( Qu ) ( \alpha ) || .
\end{array}
\eeq
In the appendix, we prove that for $B$ large enough, \[ || ( Q V P
) ( \alpha ) || \leq C_4 B^{ - 1/2 }, \] and
\[ || ( P^{*} Q ) ( \alpha ) || \leq C_5 ~| Im \alpha |~B^{ - 1/2 }
{}. \] With these estimates, we obtain from (4.29,) \beq{4.30}
\begin{array}{rl}
( a - C_6 \gamma B ^{-1} ) || ( P u ) ( \alpha ) || ^2 & \leq || (
P u ) ( \alpha ) || \: || ( P v ) ( \alpha ) || \\ & \\ & + C_7 B^{
- 1/2 } || ( P u )
( \alpha ) || \: || ( Q u ) ( \alpha ) || , \end{array}
\eeq
where the constants $ C_6 $ and $ C_7 $
depend only on $ || V ||_{ \infty } $,
$ ||\nabla V ||_{ \infty } $,
and $ || \nabla \rho ||_{ \infty } $. To estimate $ || ( Qu ) (
\alpha ) || $,
it follows from the resolvent equation and (4.27) that \beq{4.31}
\begin{array}{rl}
|| ( Qu ) ( \alpha ) || & \leq || ( Q ( H_A - z ) ^{-1} v ) (
\alpha ) || + || \{ Q ( H_A - z )^{-1} Q V ( Q + P ) u \} ( \alpha
) || \\ & \\
& \leq C_1 B ^{-1} || v ( \alpha ) || + C_1 B ^{-1} M_0 || ( Qu ) (
\alpha ) || \\ & \\ & + C_1 B ^{-1} M_0 || ( Q V P u ) ( \alpha )
|| , \end{array}
\eeq
with $ M_0 \equiv || V ||_{ \infty } < \infty $. Using the estimate
on $ Q V P $ derived in the appendix and taking $ B > 2 M_0 C_1 $,
we obtain,
\beq{4.32}
|| ( Q u ) ( \alpha ) || \leq 2 C_1 B^{ -1 } || v ( \alpha ) || +
C_8 B^{ -3/2 } || ( P u ) ( \alpha ) ||, \eeq
where $ C_8 \equiv 2 M_0 C_1 C_2 $. Substituting (4.32) into
(4.30), we obtain
\beq{4.33}
( a - C_6 \gamma B ^{-1} - C_7 C_8
B^{ -2 } ) || ( P u ) ( \alpha ) || \leq ( C_1 + 2 C_1 C_7 B^{ -
3/2 } ) || v ( \alpha )|| . \eeq
This proves (4.25) for $B$ large enough. Inserting (4.25) into
(4.32) yields (4.26). $ \Box $ \vspace{.1in}

\noindent
{\bf Corollary 4.1}. {\em Let $ {\cal O} $ be an open, connected,
bounded subset of $ I \hspace{-.17cm} R ^2 $ with smooth boundary
and suppose $ {\cal E} \subset I \hspace{-.17cm} R ^2
\setminus {\cal O} $.
Let $ E \in \sigma_0 \setminus \{ B \} $ and assume that $ ( B , E
, V ) $ satisfy (4.9) or (4.10) for some $ a > 0 $. Let $ \chi _X ,
X = {\cal O} $ and $ {\cal E} $, be bounded functions with support
in $X$
and s.t.\ $ || \chi _X ||_{ \infty }
\leq 1 $. Then,
\beq{4.34}
{\displaystyle \sup_{ \epsilon \neq 0 }} \: || \chi _{ {\cal E} } (
H_A + V - E - i \epsilon ) ^{-1} \chi _{ {\cal O} } || \leq C \max
\: \{ a ^{-1} , B ^{-1} \} e^{ - \gamma d } , \eeq
where $C$ and $ \gamma $ are as in Theorem 4.3 and $ d \equiv
\mbox{dist} \: ( {\cal O} , {\cal E} ) $.} \vspace{.1in}

\noindent
{\bf Proof}. This is an immediate consequence of Theorem 4.3. We
set $ \rho (x) \equiv \mbox{dist} \: ( x , {\cal O} ) $ and choose
$ v \equiv \chi _{ {\cal O} } v $. Then, $ e^{ i \alpha \rho } v =
v , ~~ \forall \: \alpha \in \C $. For $ u $ a solution of
$ ( H_A + V - E - i \epsilon ) u = \chi _{ {\cal O} } v $, one has
$ \forall \: \alpha \in \C , | \mbox{Im} \: \alpha | < \gamma $,
\[ \begin{array} {c}
|| \chi _{ {\cal E} } ( H_A + V - E - i \epsilon ) ^{-1} \chi _{
{\cal O} } v || =
|| \chi _{ {\cal E} } ( P + Q ) u || \\ \\ \leq e^{ -d (
\mbox{\scriptsize Im} \: \alpha ) } \{ || e^{ \alpha \rho }
P u || + || e^{ \alpha \rho } Qu || \} \\ \\ \leq e^{ - d (
\mbox{\scriptsize Im} \: \alpha ) } C \: \max \: \{ a ^{-1} ,
B ^{-1} \} || v || , \end{array} \]
by Theorem 4.3. Taking $ \mbox{Im} \: \alpha \rightarrow \gamma $,
we obtain (4.34). $ \Box $


\setcounter{chapter}{5}
\setcounter{equation}{0}

\section{Proof of the Main Theorem}\label{S.5}

We will show below that Corollary 4.1 implies
hypothesis~[H1]$(\gamma_0, \ell_0)$ of~\cite{[CoHi]}. This
hypothesis, along with Wegner's estimate, Theorem 3.1, are the main
starting points of the multi-scale analysis described
in~\cite{[CoHi]}. The goal of this analysis is to verify the main
assumption~(A2) of Theorem 3.2 of~\cite{[CoHi]}, which gives
sufficient conditions for pure point spectrum in an interval and
exponential decay of eigenfunctions. A version of Kotani's trick,
necessary to control the singular continuous spectrum for the
models studied here, follows from Lemma 3.2.

In order to make this paper more self-contained, we recall the main
points of this analysis here and refer to~\cite{[CoHi]} for the
details. To reduce the family $H_{\omega}$ in~\rf{2.1} to a one
parameter family, we consider variations $\omega' \in \Omega$ for
which only
$\lambda_0$ changes. For $\omega$ fixed and $\lambda \equiv
\lambda_0 (\omega') - \lambda_0 (\omega)$, we have \beq{5.1}
H_{\omega'} = H_{\omega} + \lambda u \equiv H_0 + \lambda u =
H_{\lambda},
\eeq

\noindent with $u$ satisfying~(V1). Let $R_{\lambda} (z) \equiv
(H_{\lambda} - z)^{-1}$. We first check the compactness condition,
(A1), of~\cite{[CoHi]}. By the diamagnetic inequality,
\beq{5.2}
e^{- H_{\lambda} (A) t} \le e^{- H_{\lambda} (0) t},\ t \ge 0 \eeq

\noindent and it is clear that $u^{1 \over 2} e^{- H_{\lambda} (0)
t} u^{1 \over 2},\ t > 0$, is compact.
Using the integral representation
\beq{5.3}
R_{\lambda} (x) = \int_0^{\infty} e^{- (H_{\lambda} (A) - x) t} d
t,\ x < 0,
\eeq

\noindent it follows from the norm-convergence of the integral and
inequality~\rf{5.2} that
$u^{1 \over 2} R_{\lambda} (z) u^{1 \over 2}$ is compact for $I m z
\not= 0,\ \forall\ \lambda$. This, for $\lambda = 0$, is
condition~(A1) of~\cite{[CoHi]}.

The second condition~(A2) is that
$\exists I_0 \subset I,\ I$ some interval, and $| I_0 | = | I |$
s.t. $\forall\ E \in I_0$,
\beq{5.4}
\build{\sup}_{\epsilon \not= 0}^{} || R_0 (E + i \epsilon) u^{1
\over 2} || < \infty .
\eeq

\noindent The multi-scale analysis is used to verify this condition
for a.e. $\omega$ (recall $H_0 = H_{\omega} $). The main theorem,
which we recall in the present context, concerning
condition~\rf{5.4}, is the following.

\begin{theorem}\label{T.5.1}

(Theorem 2.3 of~\cite{[CoHi]}). Let $\gamma_0 > 0$. $\exists$ a
minimum length scale $\ell^{\star} \equiv \ell^{\star} (\gamma_0,
C_W)$, s.t.~: if [H1]$(\gamma_0, \ell_0)$ holds at energy $E$ for
$\ell_0 > \ell^{\star}$, then for $\P$ - a.e. $\omega\ \exists$ a
finite constant $d_{\omega} > 0$ s.t.
$$
\build{\sup}_{\epsilon \not= 0}^{}
|| (H_{\omega} - E - i \epsilon)^{-1} u^{1 \over 2} || < d_{\omega}
\delta (u),
$$

\noindent where $\delta (u)$ depends only on $u$.

\end{theorem}

We prove below that [H1]$(\gamma_0, \ell_0)$ holds at each energy
in $ [ B - M_0, B - {\cal O}( B^{-1} ) ] \cup I_0 (B) \cap \sigma_0
$ with a suitable $\gamma_0$ (see Proposition 5.1) and for all
$\ell_0$ large enough provided $B$ is large. By Theorem 3.2
of~\cite{[CoHi]}, this theorem and the compactness result shown
above imply that $H_{\lambda}$ in~\rf{5.1} has pure point spectrum
in this set for a.e. $\lambda$. By the probabilistic arguments
of~\cite{[CoHi]}, we conclude that $H_{\omega}$ has only pure point
spectrum in this set for $\P$ - a.e. $\omega$.

The second main theorem which we recall here allows us to prove
exponential decay of the eigenfunctions.

\begin{theorem}\label{T.5.2}

(Theorem 2.4 of~\cite{[CoHi]}). Let $\chi_x$ be the characteristic
function of a unit cube centered at $x \in \R^2$. Under the
assumptions of Theorem~\ref{T.5.1}, for $\P$ - a.e. $\omega\
\exists$ a finite constant $d_{\omega} > 0$ s.t. for all $x$, $|| x
||$ large enough,
$$
\build{\sup}_{\epsilon > 0}^{}
|| \chi_x (H_{\omega} - E - i \epsilon)^{-1} u^{1 \over 2} || \le
d_{\omega} e^{- \gamma_1 || x ||},
$$

\noindent where $\gamma_1 \equiv (1 / 6 \sqrt 2) \gamma_0$,
$\gamma_0$ as in Theorem~\ref{T.5.1}.

\end{theorem}

Let us remark that for our problem,
$\gamma_0 \sim B^{\sigma}$, for some $\sigma > 0$ so there is
exponential decay in the $B$-field also. We now turn to the proof
of [H1]$(\gamma_0, \ell_0)$.

To begin, we introduce some geometry. In this section, we work with
subregions of the lattice $\Gamma \equiv e^{i \pi / 4} \sqrt 2
\Z^2$, introduced in section~\ref{S.4}, rather than in $\Z^2$.
Recall that there is a 1:1 correspondence between bonds $b_j \in
\Gamma$ and vertices $j \in \Z^2$. We arbitrarily choose a vertex
of $\Gamma$ as the origin and define boxes $\Lambda_{\ell} \subset
\Gamma$ relative to this point,
$$
\Lambda_{\ell} \equiv \left\{ x \in \R^2 | | x_{i} | \le \ell / 2
{}~~for ~~i = 1, 2 \right\}.
$$

\noindent For convenience, we fix points so the bond $b_0$ has one
of its ends at $0 \in \Gamma$. For any $\delta > 0$, consider
$\Lambda_{\ell, \delta} \equiv \left\{ x \in \Lambda_{\ell} | dist
(x, \partial \Lambda_{\ell}) < \delta \right\}$. Let $\chi_{\ell,
\delta}$ be the $C^2$-function which satisfies $\chi_{\ell, \delta}
> 0$, $| \nabla \chi_{\ell, \delta} | \subset \Lambda_{\ell}
\backslash \Lambda_{\ell, \delta}$ and
$\chi_{\ell, \delta} | \Lambda_{\ell, \delta} = 1$. Let $W (\chi)
\equiv [\chi, H_A]$, for any $\chi \in C^2$. Let $V_{\Lambda}
\equiv V | \Lambda$, $\Lambda \subset \R^2$ and $H_{\Lambda} \equiv
H_A + V_{\Lambda}$, as in section 4.

We apply the multi-scale analysis to $H_{\Lambda}$ relative to the
lattice $\Gamma$. We verify condition~[H1]$(\gamma_0, \ell_0)$
of~\cite{[CoHi]} using Corollary 4.1 and the geometric resolvent
equation ($G R E$). We must show that for $E \in [ B - M_0, B -
{\cal O}( B^{-1}) ] \cup I_0 (B) \cap \sigma_0 $ and for all
$\ell_0$ sufficiently large, that the following holds:

\bigskip

[H1] $( \gamma_0, \ell_0 )$ For some $\gamma_0 > 0,\ \ell_0 > 1,
{}~~\exists \xi > 4 $ s.t.
$$
\P \left\{ \build{\sup}_{\epsilon > 0}^{} || W (\chi_{\ell,
\delta}) R_{ \Lambda_{\ell_0} } (E + i \epsilon) \chi_{ \ell_{0 /
3} } || \le e^{- \gamma_0 \ell_0} \right\} \ge 1 - \ell_0^{- \xi}.
$$

We begin with a simple lemma which allows us to control the
gradient term in $W (\chi_{\ell, \delta})$.

\begin{lemma}\label{L.5.3}

Let $H_{\Lambda} \equiv (p - A)^2 + V_{\Lambda}$ and write
$R_{\Lambda} \equiv R_{\Lambda} (E + i \epsilon),\ \epsilon \not=
0,\ E \in \R $. For any $u \in L^2 (\R^2),\ || u || = 1$, we have
for $i = 1,2$,
\beq{5.5}
|| (p - A)_i R_{\Lambda} u ||^2 \le || R_{\Lambda} u || + (2 M_0 +
| E |) || R_{\Lambda} u ||^2,
\eeq

\noindent where $M_0 \equiv || V_{\Lambda} ||_{\infty} > 0$.
Moreover, for any bounded $\chi \in C^1$, we have, \beq{5.6}
\ba{ll}
\Sum_{i = 1}^2 || \chi (p - A)_i R_{\Lambda} u ||^2 \le || \chi
R_{\Lambda} u || & + (2 M_0 + | E |) || \chi R_{\Lambda} u ||^2
\\[3mm]
& + 2 \Sum_{i = 1}^2 || (\partial_i \chi) u ||~~|| \chi (p - A)_i
R_{\Lambda} u ||.
\ea
\eeq

\end{lemma}

\begin{prf}

The inequality \rf{5.4} follows directly from the equality $$
\langle R_{\Lambda} u,\ H_A R_{\Lambda} u \rangle = \langle
R_{\Lambda} u,\ u \rangle - \langle R_{\Lambda} u,\ (V_{\Lambda} -
E - i \epsilon) R_{\Lambda} u \rangle ,
$$

\noindent and the Cauchy-Schwartz inequality. The
inequality~\rf{5.5} follows in the same way by writing out $|| \chi
(p - A)_i R_{\Lambda} u ||^2$. \end{prf}

We now prove the main result of this section. Recall that $ \cal R
$ denotes the ribbon defined in section 4.

\begin{proposition}\label{P.5.4}

Let $\chi_2$ be any function, $|| \chi_2 ||_{\infty} \leq 1$,
supported on $ \Lambda_{\ell} \cap Ext {\cal R} $, where $ ~~Ext
{\cal R} \equiv \left\{ x \in \R^2 | \lambda x \not\in {\cal R
}~\forall\ \lambda \ge 1 \right\}$, so that, in particular, $\supp
\chi_2 \cap {\cal R } = \emptyset$. For any $E \in \sigma_0
\backslash \{ B \},\ \delta > 0, \ \epsilon > 0,$ and $a > 0$, we
have \beq{5.7}
\ba{ll}
\build{\sup}_{\epsilon \not= 0}^{} || \chi_2 R_{\Lambda_{\ell}} (E
+ i \epsilon) \chi_{\ell / 3} || \le & C e^{- \gamma d} \max
\left\{ a^{-1}, B^{-1} \right\} \cdot \max \left\{ \delta^{-1},
\right. \\[3mm]
& \left. (2 M_0 + | E |) \delta^{-2} \right\}, \ea
\eeq

\noindent where $C$ and $\gamma$ are as in Theorem 4.3 and $d
\equiv r_1 - 3 \epsilon ~~~(r_1 \equiv diam \cal R)$, with a
probability larger than
\beq{5.8}
1 - \left\{ C e^{- m \ell} + C_W [ \dist (E,B) ]^{-2} || g
||_{\infty} \delta B \ell^2 \right\}.
\eeq

\noindent In particular, for $\chi_{\ell, \delta}$ defined above
and $E \in \sigma_0$ with $a = {E - B \over 2} = {\cal O} \left(
B^{- 1 + \sigma}
\right)$, any $\sigma > 0$, we have that for any $\ell_0 > \sqrt 2$
and large enough, and any\break
$\xi > 4, ~~\exists ~~B (\ell_0) > 0$ s.t. $\forall\ B > B
(\ell_0)$, [H1] $(\gamma_0,\ \ell_0)$ holds for some $\gamma_0 >
\gamma d / 4 \ell_{0} > 0$, so that $ \gamma_{0} = {\cal O}\{ min (
B^{1/2}, B^{\sigma} ) \}. $

\end{proposition}

\begin{prf}

{\bf 1.} By Corollary 4.1, $\exists ~~B_0$ s.t. $B > B_0$ implies
$\exists$ a ribbon ${\cal R} \subset \Lambda_{\ell} \backslash
\Lambda_{\ell / 3}$ (with a probability given by~\rf{4.7})
satisfying
\beq{5.9}
\dist({\cal R }, \partial \Lambda_{\ell}), ~~and ~~dist( {\cal R},
\partial \Lambda_{\ell /
3}) > {1 \over \sqrt 2} + r_u > 0,
\eeq
and
\beq{5.10}
r_1 \equiv diam {\cal R} > 2 \left( {1 \over \sqrt 2} - r_u
\right), \eeq
and such that
\beq{5.11}
V (x) + B - E > -a \ ~~\forall\ x \in {\cal R},\ a = {E - B \over
2} \cdotp \eeq

\noindent (We assume $E > B$ ; similar arguments hold for $E < B$).
For any\break
$\epsilon > 0,\ 3 \epsilon \ll r_1$, define the border of $\cal R$
by $$
{\cal R}_{\epsilon} \equiv \left\{ x \in {\cal R} | ~~\dist(x,
\partial {\cal R}) <
\epsilon \right\}.
$$

\noindent Then ${\cal R}_{\epsilon} \equiv {\cal R}_{\epsilon}^+
\cup {\cal R}_{\epsilon}^-$, where
${\cal R}_{\epsilon}^{\pm}$ are two disjoint, connected subsets of
$\cal R $. Let
${\cal C} \equiv \left\{ x \in {\cal R} | \dist (x, {\cal
R}_{\epsilon}^+) = \dist (x,
{\cal R}_{\epsilon}^-) \right\}$ ; $\cal C$ is a closed, connected
path in $\cal R$. Let ${\cal C}_{\epsilon} \equiv \{ x \in {\cal R}
| \dist (x, {\cal C} ) < \epsilon \} \subset {\cal R}$ and
\beq{5.12}
\dist \left( {\cal C}_{\epsilon}, {\cal R}_{\epsilon}^{\pm} \right)
\ge r_1 - 3 \epsilon.
\eeq

\noindent This is strictly positive. Because of this, we can adjust
${\cal C}_{\epsilon}$ so that $\partial {\cal C}_{\epsilon}$ is
smooth. We need two, $ C^2$, positive cut-off functions. Let
$\chi_{\cal R} > 0$ sa\-tis\-fy $\chi_{\cal R} | {\cal
C}_{\epsilon} = 1$ and $\supp | \nabla \chi_{\cal R} | \subset
{\cal R}_{\epsilon}$. Let $\chi_1$ satisfy
$\chi_1 | \Lambda_{\ell / 3} = 1$ and $\supp | \nabla \chi_1 |
\subset {\cal C}_{\epsilon}$ (see Figure 4). By simple commutation,
we have (with
$\chi_2$ as in the proposition),
\beq{5.13}
\ba{ll}
\chi_2 R_{\Lambda_{\ell}} (E + i \epsilon) \chi_{\ell / 3} & =
\chi_2 R_{\Lambda_{\ell}} \chi_1 \chi_{\ell / 3} \\[3mm] & = \chi_2
R_{\Lambda_{\ell}} W (\chi_1) R_{\Lambda_{\ell}} \chi_{\ell / 3}
\\[3mm] & = \chi_2 R_{\Lambda_{\ell}}
\chi_{\cal R} W (\chi_1) R_{\Lambda_{\ell}} \chi_{\ell / 3}. \ea
\eeq

\noindent Next, denote by $R_{\cal R}$ the resolvent of $H_{\cal
R}$ defined
in section 4.2. The $G R E$ relating $R_{\Lambda_{\ell}}$ and
$R_{\cal R}$ is
\beq{5.14}
R_{\Lambda_{\ell}} \chi_{\cal R} = \chi_{\cal R} R_{\cal R} +
R_{\Lambda_{\ell}} W
(\chi_{\cal R}) R_{\cal R}.
\eeq

\noindent Substituting~\rf{5.13} into~\rf{5.12} and noting that
$\chi_2 \chi_{\cal R} = 0$, we obtain
\beq{5.15}
\chi_2 R_{\Lambda_{\ell}} \chi_{\ell / 3} = \chi_2
R_{\Lambda_{\ell}} W (\chi_{\cal R}) R_{\cal R} W (\chi_1)
R_{\Lambda_{\ell}} \chi_{\ell / 3}.
\eeq

\noindent Note that from~\rf{5.11} and the choice of $\chi_{\cal
R}$ and $ \chi_1$, we obtain that
\beq{5.16}
\dist \left( \supp W (\chi_{\cal R}),\ \supp W (\chi_1) \right) \ge
r_1 - 3 \epsilon.
\eeq

\noindent We apply Wegner's estimate, Theorem~\ref{T.3.1}, to
control the two
$R_{\Lambda_{\ell}}$ factors in~\rf{5.14}, and the decay estimate,
Corollary 4.1, to control the factor $R_{\cal R}$, which is
possible due to the localization of $W (\chi_{\cal R})$ and $W
(\chi_1)$ and~\rf{5.16}.

\bigskip

{\bf 2.} To estimate the $R_{\cal R} (E + i \epsilon)$
contribution, we use Corollary 4.1 with ${\cal O} \equiv {\cal
C}_{\epsilon}$ and ${\cal E } = {\cal R}_{\epsilon}$. Let $\chi_X,
{}~~X = {\cal O}$ and ${\cal E}$, be a characteristic function on
these sets. Then
$W (\chi_{\cal R})
\chi_{\cal E} = W (\chi_{\cal R})$ and $\chi_{\cal O} W (\chi_1) =
W (\chi_1)$. Inserting these localization functions into~\rf{5.15},
we obtain from Corollary 4.1,
\beq{5.17}
|| \chi_{\cal E} R_{\cal R} (E + i \epsilon) \chi_{\cal O} || \le C
\max
\left\{ a^{-1}, B^{-1} \right\} e^{- \gamma d}, \eeq

\noindent with probability larger than
\beq{5.18}
1 - C e^{- m \ell},
\eeq

\noindent for some $m = m (1 - p) > 0$ and $0 < C < \infty$. The
factor $d$ satisfies
\beq{5.19}
d \ge r_1 - 3 \epsilon,
\eeq

\noindent where $r_1 \equiv ~diam {\cal R} $ as in~\rf{5.10}.

\bigskip

{\bf 3.} Next, we turn to
\beq{5.20}
W (\chi_1) R_{\ell} (E + i \epsilon) \chi_{\ell / 3}, \eeq
and
\beq{5.21}
\chi_2 R_{\ell} (E + i \epsilon) W (\chi_{\cal R}), \eeq

\noindent where we write $ R_{\ell} $ for $ R_{ \Lambda_{\ell}}$
for short. We will bound $R_{\ell}$ by Wegner's estimate and the
terms~\rf{5.20}-\rf{5.21} via Lemma~\ref{L.5.3}. From
Theorem~\ref{T.3.1}, we have for any $\delta > 0$, \beq{5.22}
|| R_{\ell} (E + i \epsilon) || < \delta^{-1}, \eeq

\noindent with probability larger than
\beq{5.23}
1 - C_W [ \dist (E, B) ]^{-2} || g ||_{\infty} \delta B \ell^2.
\eeq

\noindent From~\rf{5.22} and Lemma 5.1, both~\rf{5.20} and
\rf{5.21} are bounded above by
\beq{5.24}
2^{1 \over 2} \max \left\{ \delta^{- 1 \over 2},\ ( 2 M + | E |
)^{1 \over 2} \delta^{-1} \right\},
\eeq

\noindent with probability at least~\rf{5.23}.

\bigskip

{\bf 4.} Using the estimate $P (A \cap B) \ge P (A) + P (B) -1$,
and \rf{5.17}-\rf{5.18} and \rf{5.23}-\rf{5.24}, we find \beq{5.25}
\ba{ll}
|| \chi_2 R_{\ell} (E + i \epsilon) \chi_{\ell / 3} || \le & 2 C
\max \left\{ a^{-1}, B^{-1} \right\} \cdot \max \left\{
\delta^{-1}, \right. \\[3mm]
& \left. (2 M_0 + | E |) \delta^{-2} \right\} \cdot e^{- \gamma d},
\ea
\eeq

\noindent with probability at least
\beq{5.26}
1 - \left\{ C e^{- m \ell} + C_W [ \dist (E, B) ]^{-2} || g
||_{\infty} \delta B \ell^2 \right\}.
\eeq

\noindent This proves the first part of the proposition.

\bigskip

{\bf 5.} To estimate $W \left( {\chi}_{\ell, \delta} \right)
R_{\ell} \chi_{\ell / 3}$, we use the second formula of
Lemma~\ref{L.5.3}, \rf{5.6}, which gives \beq{5.27}
\ba{ll}
|| \chi_2 (p - A)_i R_{\ell} \chi_{\ell / 3} ||^2 & \le || \chi_2
R_{\ell} \chi_{\ell / 3} || + (2 M + | E |) || \chi_2 R_{\ell}
\chi_{\ell / 3} ||^2 \\[3mm] & + 2 \build{\max}_{i = 1,2}^{} ||
(\partial_i \chi_2) R_{\ell}
\chi_{\ell / 3} || ~~|| \chi_2 (p - A)_i R_{\ell} \chi_{\ell /
3}||. \ea
\eeq

\noindent Since $\partial_i \chi_2$ satisfies the same condition as
$\chi_2$, the factor $|| (\partial_i \chi_2) R_{\ell} \chi_{\ell /
3} ||$ in~\rf{5.27} satisfies the estimate~\rf{5.24} with possibly
a different constant. Sol\-ving the quadratic inequality~\rf{5.27},
we obtain
\beq{5.28}
\ba{ll}
|| \chi_2 (p - A)_i R_{\ell} \chi_{\ell / 3} || \le &
\build{\max}_{i = 1,2}^{}
\left\{ \vphantom{ \left. \left. \chi_{\ell / 3} ||^2 \right)
\right]^{1 / 2} }
|| (\partial_i \chi_2) R_{\ell} \chi_{\ell / 3}

|| \right. \\[3mm]
& + \left[ || (\partial_i \chi_2) R_{\ell} \chi_{\ell / 3} ||^2 +
\left( || \chi_2 R_{\ell} \chi_{\ell / 3} || \right. \right.
\\[3mm]
& \left. \left. \left. \! \! + (2 M_0 + | E |) || \chi_2 R_{\ell}
\chi_{\ell / 3} ||^2 \right) \right]^{1 / 2} \right\},
\ea
\eeq

\noindent which can be estimated as in~\rf{5.25}. Finally, we write
\beq{5.29}
|| W \left( {\chi}_{\ell, \delta} \right) R_{\ell} \chi_{\ell / 3}
|| \le || \left( \Delta {\chi}_{\ell, \delta} \right) R_{\ell}
\chi_{\ell / 3} || + 2 \Sum_{j = 1}^{2} || \left( \partial_j
{\chi}_{\ell,\delta} \right) (p - A)_j R_{\ell} \chi_{\ell / 3} ||,
\eeq

\noindent which can be estimated from~\rf{5.25} with $\chi_2 \equiv
\Delta {\chi}_{\ell, \delta}$ and~\rf{5.27} with $\chi_2 \equiv
\left( \partial_j {\chi}_{\ell, \delta} \right)$.

\bigskip

{\bf 6.} We now show that for any $\ell_{0} $ large enough, $
\exists B_{0} \equiv B_{0}( \ell_{0} ) $ such that for all $ B >
B_{0} $,
condition $ [H1]( \ell_{0}, \gamma_{0}) $ is satisfied with $
\gamma_{0} = {\cal O}\{ min (B^{1/2}, B^{\sigma} )\} $. We take $E
\in [ B - M_0 , B - {\cal O}( B^{-1} ) ] \cup I_0 (B) \cap
\sigma_{0} $ and $a = {E - B \over 2} = {\cal O} \left( B^{- 1 +
\sigma} \right)$, for any $\sigma > 0$. First, we require that
(5.25) be bounded above by $e^{- \gamma d / 2} $. This leads to the
condition \beq{5.30}
C \delta^{-2} B^{ 2 - \sigma} e^{ - \gamma d } \leq e^{- \gamma d /
2},
\eeq
where $\gamma = C_{2} min \{ B^{1 \over 2}, B^{\sigma} \}$. This
condition implies that we must choose $ \delta $ in (5.22) to
satisfy
\beq{5.31}
\delta > B^{ 1 - ( \sigma/2)} e^{ - \gamma d / 4}. \eeq
If we now define
$$
\gamma_{0} \equiv \gamma d / 4 \ell_{0}
$$
we find that
$$
|| W \left( {\chi}_{\ell, \delta} \right) R_{\ell} \chi_{\ell / 3}
|| \le e^{ - \gamma_{0} \ell_{0} }. $$

\noindent
Next, the probability estimate (5.26) leads to the condition
\beq{5.32}
C e^{ - m \ell_{0} } + C_{2} B^{3 - 2 \sigma } \delta \ell^{2}_{0}
\leq \ell_{0}^{ - \xi },
\eeq
or, for all $ \ell_{0} $ large,
\beq{5.33}
C_{3} B^{3 - 2 \sigma } \delta \ell^{2}_{0} \leq \ell_{0}^{ - \xi
},
\eeq
for some $ \xi > 4 $.
We can choose $ \delta $ so that both conditions (5.31) and (5.33)
are satisfied provided the condition
\beq{5.34}
\ell_{0}^{ \xi + 2 } < B^{3/2 - (5/2) \sigma } e^{ \gamma d / 4},
\eeq
is satisfied for some $ \xi > 4 $. It is clear from the definition
of $ \gamma $, that for any $ \ell_{0} $, there exists a $B_{0}
\equiv B_{0}( \ell_{0} ) $ such that condition (5.34) is satisfied
for all $ B > B_{0} $. This completes the proof of the theorem.
\end{prf}

\newpage


\setcounter{chapter}{6}
\setcounter{equation}{0}

\section{Appendix}

\noindent
The following estimates hold for all $ B $ sufficiently large.
\vspace{.1in}

\noindent
{\bf Lemma A.1}. {\em Let $ V \in C_b^2 ( I \hspace{-.17cm} R ) $.
$ \exists $ constant $ C > 0 $ depending only on $ || \partial^{
\alpha } V ||_{ \infty }, ~~| \alpha | = 0,1,2 $, such that
$ \forall \: \alpha \in S $},
\beq{A.1}
|| P ( \alpha ) V Q ( \alpha ) || \leq C B^{ -1/2 } . \eeq
\vspace{.1in}

\noindent
{\bf Proof}. Let $ z \equiv B - 1 $, so $ z \in \rho ( H_A ( \alpha
)) $ for $ \alpha \in S $. We have \beq{A.2}
\begin{array}{rl}
P ( \alpha ) V Q ( \alpha ) & = P ( \alpha ) ( H_A ( \alpha ) - z )
( H_A ( \alpha ) - z ) ^{-1} V Q ( \alpha ) \\ & \\
& = P ( \alpha ) ( H_A ( \alpha ) - z ) V ( H_A ( \alpha ) - z )
^{-1}
Q ( \alpha ) \\ & \\
& + P ( \alpha ) ( H_A ( \alpha ) - z ) ( H_A ( \alpha ) - z )
^{-1}
[ V , H_A ( \alpha ) ] ( H_A ( \alpha ) - z ) ^{-1} Q ( \alpha ).
\end{array}
\eeq
Recall that $ P ( \alpha ) $ is analytic in $ \alpha \in S $. As
\beq{A.3}
\begin{array}{rl}
P ( \alpha ) ( H_A ( \alpha ) - z ) & = P ( \alpha ) ( B - z ) \\ &
\\
& = P ( \alpha ) , \end{array}
\eeq
for $ \alpha \in I \hspace{-.17cm} R $, the identity principle for
analytic funct
ions
implies this holds for $ \alpha \in S $. This result (A.3) and
estimates (4.16)-(4.17) imply that the first term on the right in
(A.2) is bounded as
\[ || P ( \alpha ) ( H_A ( \alpha ) - z ) V ( H_A ( \alpha ) - z )
^{-1}
Q ( \alpha ) || \leq C_0 || V ||_{ \infty } B ^{-1} , \] $ \forall
\: \alpha \in S $. As for the second term, the commutator
is (see (4.18))
\[ [ V , H_A ( \alpha ) ] = 2 i ( p - \alpha \nabla \rho - A )
\cdot
\nabla V - \Delta V . \]
The resulting term in (A.2) involving $ \Delta V $ is treated as
above. As for the derivative term, it suffices to show \beq{A.4}
|| ( H_A ( \alpha ) - z ) ^{-1} ( p_i - \alpha \partial _i \rho -
A_i ) || \leq C_1 B^{ 1/2 } , \eeq
for all $ \alpha \in S $. To see this, let $ V_i ( \alpha ) \equiv
( p_i - \alpha
\partial _i \rho - A_i ) $
and $ R ( \alpha ) \equiv ( H_A ( \alpha ) - z ) ^{-1} $. For $ u =
R ( \alpha ) v , ~~ v \in L^2
( I \hspace{-.17cm} R ^2 ) $,
we have
\[ \sum_{ i = 1 }^2 \: || V_i ( \alpha ) u ||^2 = \langle u , \{
H_A ( \alpha ) +
2i ( \mbox{Im} \: \alpha ) \nabla \rho \cdot V ( \alpha ) \}
u \rangle \]
\[ = \langle R ( \alpha ) v , v \rangle + z || u ||^2 + 2i
( \mbox{Im} \: \alpha ) \langle u , \nabla \rho \cdot V ( \alpha )
u \rangle . \]
This leads to a quadratic inequality for each $ i = 1 , 2 $,
\[ \begin{array}{rl}
|| V_i ( \alpha ) R ( \alpha ) v ||^2 & \leq || R ( \alpha ) v || (
|| v || + | z
| \:
|| R ( \alpha ) v || ) \\ & \\
& ~~~~ + 2 | \mbox{Im} \: \alpha | \: || \nabla \rho ||_{ \infty }
|| R ( \alpha ) v || \left \{ {\displaystyle \max_{ j = 1,2 }} \:
|| V_j ( \alpha ) \right. \\ & \\
& \left. ~~~~ \times R ( \alpha ) v || \right \} . \end{array} \]
Solving this, and noting that $ | \mbox{Im} \: \alpha | \leq B^{
1/2 } , ~~
|z| = {\cal O} (B) $,
and, for this $ z , ~~ || R ( \alpha ) || < C_0 $ by (4.16)-(4.17),
we get
\beq{A.5}
|| V_i ( \alpha ) R ( \alpha ) || \leq C_1 B^{ 1/2 } , \eeq
which is (A.4). $ \Box $
\vspace{.1in}

\noindent
{\bf Lemma A.2}. {\em Let $\rho $ be the distance function defined
in section 4.
$ \exists $ constant $ C > 0 $ depending only on $ || \partial^{
\alpha } \rho ||_{ \infty }, ~~| \alpha | = 0,1,2 $, such that
$ \forall \: \alpha \in S $},
\beq{6.6}
|| P ( \alpha )^{\star} Q ( \alpha ) || \leq C B^{ -1/2 }|
Im{\alpha}|,
\eeq
\noindent
for $ | Im {\alpha} | \leq B^{ 1/2 }$.
\vspace{.1in}

\noindent
{\bf Proof.} We can assume that $ \alpha $ is purely imaginary by a
standard unitary equivalence argument (note that when $\alpha$ is
real, $ || P ( \alpha )^{\star} Q ( \alpha ) || = 0 $.) Let $ z
\equiv B - 1 $, as in Lemma A.1. We then have by (6.3), \beq{6.7}
\begin{array}{rl}
P ( \alpha )^{\star} Q ( \alpha ) = & P ( \alpha )^{\star} (
H_A^{\star} ( \alpha ) - z )^{-1} Q ( \alpha ) \\ & \\ = & P (
\alpha )^{\star}\{ ( H_A^{\star}( \alpha ) - z )^{-1} -
( H_A ( \alpha ) - z )^{-1} \} Q ( \alpha ) \\ & \\ & + P ( \alpha
)^{\star}( H_A ( \alpha ) - z )^{-1} Q ( \alpha ). \end{array}
\eeq
The last term is $ {\cal O }( B^{-1} ) $ by (4.16)-- (4.17). Let us
write $ \alpha \equiv i \beta $, with $ \beta $ real. Then by the
resolvent equation, we have
\beq{6.8}
(H_A^{\star}( \alpha ) - z )^{-1} - ( H_A ( \alpha ) - z )^{-1} =
-2i{\beta}( H_A^{\star}( \alpha ) - z )^{-1} [ 2 ( p - A ) \cdot
{\nabla } \rho + i \Delta \rho ] ( H_A ( \alpha ) - z )^{-1}. \eeq
Using (6.8) in (6.7), we obtain the bound, \beq{6.9}
\begin{array}{rl}
|| P ( \alpha )^{\star} Q ( \alpha ) || \leq & 2| Im(\alpha) | ~||
\{ 2 ( p - A ) \cdot {\nabla } \rho + i \Delta \rho \}
P^{\star}(\alpha) ||~|| ( H_A ( \alpha ) - z )^{-1} \} Q ( \alpha )
\\
& \\
\leq & C_{1}| Im(\alpha) | B^{-1} ( || 2 ( p - A ) \cdot
\nabla \rho P^{\star}(\alpha) || + C_{2} ), \end{array}
\eeq
where the constants are bounded as in the lemma. We again used
(4.16)--(4.17) and the boundedness of $ \Delta \rho $. The proof
will follow from (6.8) once we show that \beq{6.10}
|| | p - A | P^{\star}(\alpha) || \leq C_{1} B^{ 1/2 }. \eeq
Inequality (6.10) follows directly as in (6.4) if we write
\beq{6.11}
( H_A^{\star} ( \alpha ) - z )^{-1} ( p_i - A_{i} ) = ( H_A^{\star}
( \alpha ) - z ) ^{-1}\{ ( p_i + i \beta \partial _i \rho - A_i ) -
i \beta \partial_i \rho \}, \eeq
and note that $ \beta \leq B^{ 1/2 } $. This proves (6.10). $ \Box
$



\end{document}